\pdfoutput=1
\RequirePackage{ifpdf}
\ifpdf 
\documentclass[pdftex]{sigma}
\else
\documentclass{sigma}
\fi

\numberwithin{equation}{section}

\newtheorem{Theorem}{Theorem}[section]
\newtheorem*{Theorem*}{Theorem}

\newtheorem{Lemma}[Theorem]{Lemma}
\newtheorem{Proposition}[Theorem]{Proposition}

\theoremstyle{definition}

\newtheorem{Remark}[Theorem]{Remark}

\newcommand{\Z}{\mathbb{Z}}
\newcommand{\R}{\mathbb{R}}
\newcommand{\C}{\mathbb{C}}
\newcommand{\K}{\mathbb{K}}
\newcommand{\QQ}{\mathbb{Q}}
\newcommand{\g}{\mathfrak{g}}
\newcommand{\h}{\mathfrak{h}}
\newcommand{\res}{\operatorname{res}}
\newcommand\Wam{W^{(a)}_{m}}
\newcommand{\rwma}{\res \Wam}
\newcommand\Qam{Q^{(a)}_{m}}
\newcommand\Fam{F^{(a)}_{m}}
\newcommand{\Q}[2]{Q^{(#1)}_{#2}}
\newcommand{\PP}[2]{P^{(#1)}_{#2}}
\newcommand\Qa{\mathcal{Q}^{(a)}}
\newcommand\Pa{\mathcal{P}^{(a)}}
\newcommand\Pam{P^{(a)}_{m}}
\newcommand\uqg{U_q(\g)}
\newcommand\uqghat{U_q(\widehat{\g})}
\newcommand\posR{\Delta^{+}}
\newcommand\La{\Lambda}
\newcommand\la{\lambda}

\usepackage[mathscr]{eucal}
\usepackage{enumitem}

\begin{document}

\allowdisplaybreaks

\newcommand{\arXivNumber}{1901.00104}

\renewcommand{\PaperNumber}{064}

\FirstPageHeading

\ShortArticleName{On Polyhedral Formulas for Kirillov--Reshetikhin Modules}

\ArticleName{On Polyhedral Formulas\\ for Kirillov--Reshetikhin Modules}

\Author{Chul-hee LEE}

\AuthorNameForHeading{C.-h.~Lee}

\Address{June E Huh Center for Mathematical Challenges, Korea Institute for Advanced Study,\\
 Seoul 02455, South Korea}
\Email{\mail{chlee@kias.re.kr}}

\ArticleDates{Received December 24, 2025, in final form June 26, 2026; Published online July 05, 2026}

\Abstract{We propose a~method to prove a~polyhedral branching formula for Kirillov--Reshetikhin (KR) modules over an untwisted quantum affine algebra. When the underlying simple Lie algebra is of exceptional type, such a~formula remains conjectural in many cases. Using a~linear recurrence relation satisfied by the characters of KR modules, we convert the verification of a~polyhedral formula into an identity between two rational functions of a~single variable with only simple poles at known locations. It is then sufficient to compare the residues at those poles, which are explicitly computable quantities. By applying this strategy, we obtain new, computer-assisted and easily verifiable proofs of known polyhedral formulas in types $F_4$ and $G_2$ within a~uniform framework.}


\Classification{17B37; 81R10}

\section{Introduction}

Let $\g$ be a~complex simple Lie algebra.
The Kirillov--Reshetikhin (KR) modules constitute an important family of finite-dimensional irreducible representations of the quantum affine algebra~$\uqghat$.
For each node $a$ of the Dynkin diagram of $\g$ and $m\in \Z_{\geq 0}$, we use the notation \smash{$W^{(a)}_m$} for the corresponding Kirillov--Reshetikhin module and \smash{$\res W^{(a)}_m$} for its restriction to a~$\uqg$-module.
An interesting problem is to understand how such a~module, or a~tensor product of such modules, decomposes into irreducible $U_q(\mathfrak{g})$-modules.
Kirillov and Reshetikhin introduced the fermionic formula~\cite{MR906858} to describe such a~decomposition; it was originally formulated as a~conjecture based on insights from the Bethe ansatz in quantum spin chains, and is often interpreted as a~combinatorial completeness of the Bethe ansatz.
This conjecture has since been proved through a~series of works including, among others, \cite{MR2428305,MR1745263, MR2254805, MR1993360}, and it provides an answer to this question by expressing the multiplicity of each irreducible summand using a~certain combinatorial rule.
Kleber's algorithm gives an efficient way to compute branching rules from the fermionic formula for fixed inputs~\cite{MR1436775}; nevertheless, for a~family of single KR modules depending on the parameter $m$, it is still desirable to have a~closed decomposition formula that is more explicit and computationally cheaper to use.
In this sense, Kleber's algorithm is complementary to the present approach: computations based on the fermionic formula, including Kleber's algorithm, can provide useful data for guessing closed formulas, while the method of this paper is designed to verify a~proposed closed formula uniformly in the parameter $m$.

If $\g$ is of classical type, there is a~well-known explicit formula often called the domino removal rule, which was proved in~\cite{MR1836791}. It is a~polyhedral formula in the sense that the highest weight of an irreducible summand with non-zero multiplicity is characterized as a~lattice point in a~suitable bounded polyhedron.
Even when $\g$ is of exceptional type, a~polyhedral formula still seems to exist, although an irreducible summand with a~multiplicity greater than one may appear.
Conjectural polyhedral formulas with multiplicities were proposed in~\cite{MR1745263} for several exceptional nodes, and many of these conjectures have since been verified, notably in~\cite{MR4114424}.
However, explicit polyhedral formulas are still not available for all exceptional nodes.
Throughout the paper, we use the same convention for enumerating the nodes of the Dynkin diagram as in~\cite{MR1745263}.
The known explicit polyhedral formulas with nontrivial multiplicities are those for type $G_2$ \cite{MR2238884}, and the cases~${(F_4,a=2)}$, $(E_6,a=3)$, $(E_7,a=2)$, and $(E_8,a=2)$ treated in~\cite{MR4114424}.

In this paper, we propose a~general method to prove a~polyhedral formula. The key objects in our approach are the coefficients that appear when the characters of KR modules are written in some exponential form.
They are essentially the residues at the poles of the generating function of the characters of KR modules.
It turns out that it is possible to decompose a~polyhedral formula into a~finite list of identities involving these coefficients; see~\eqref{eqn:CPCQ}.
Our method is well suited for a~computer-aided mechanical approach.
While the difficulty of the actual implementation varies depending on the type of $\g$, such an approach is well justified in view of our focus on exceptional types.
The method should be viewed as a~verification procedure for proposed polyhedral formulas.
Finding a~polyhedral formula remains a~separate problem: one still has to express the multiplicities as step-polynomials on suitable polyhedral regions.
The method in the present paper is formulated in the untwisted affine setting, since it uses the untwisted $Q$-system and the linear recurrence relation from~\cite{MR3948941}.
We expect the same strategy to apply to twisted types as well, but in this paper we restrict ourselves to untwisted quantum affine algebras.

After presenting the general strategy, we consider special cases in which $\g$ is of type $F_4$ or $G_2$.
Although the results stated below are already known in the literature, we obtain a~computer-aided and easily verifiable proof of the following two polyhedral formulas conjectured in~\cite{MR1745263} using the method introduced in this paper.
\begin{Theorem}[\cite{MR4114424}]\label{thm:main}
Let $\g$ be of type $F_4$. For every $m\in \Z_{\geq 0}$, the following holds:
\begin{equation}\label{eqn:sumef}
\res W^{(2)}_m =
\bigoplus_{
	\substack{
		j_1 + 2 j_2 + j_3 + j_4 \le m\\
		j_1,j_2,j_3,j_4 \in \Z_{\geq 0}
	}}
 p(j_1,j_2,j_3,j_4;m) L\left(j_1\la_1 +j_2\la_2 +j_3 \la_3 + j_4 \la_4 \right),
\end{equation}
where
\[
p(j_1,j_2,j_3,j_4;m)=\min \left( 1+j_3,1+m-j_1-2 j_2-j_3-j_4 \right) (j_4+1)
\]
and
$(\la_1,\la_2,\la_3,\la_4)=(2\omega_4,2\omega_3,\omega_2,\omega_1)$.
\end{Theorem}

\begin{Theorem}[\cite{MR2238884}]\label{thm:main2}
Let $\g$ be of type $G_2$. For every $m\in \Z_{\geq 0}$, the following holds:
\begin{equation}\label{eqn:sumeg}
\res W^{(2)}_m =
\bigoplus_{
	\substack{
		j_0+3j_1 + j_2 = m\\
		j_0,j_1,j_2 \in \Z_{\geq 0}
	}}
 p(j_0,j_1,j_2;m) L (j_1\omega_1 +j_2\omega_2 ),
\end{equation}
where
\begin{equation}\label{eqn:step-poly}
p(j_0,j_1,j_2;m)=\left\{
\min\left(
	\left\lfloor \frac{j_0}{3} \right\rfloor,
	\left\lfloor \frac{j_2-j_0}{3} \right\rfloor + \left\lfloor \frac{j_0}{3} \right\rfloor
	\right)
+1 \right\}( j_1+1 ).
\end{equation}
\end{Theorem}

This paper is organized as follows.
In Section~\ref{sec:background}, we review the necessary background for our approach such as the $Q$-system, and linear recurrence relations satisfied by the characters of KR modules.
In Section~\ref{sec:framework}, we explain the steps for proving a~polyhedral formula for KR modules.
In Section~\ref{sec:f4}, we follow the procedures described in Section~\ref{sec:framework} to give a~proof of Theorem~\ref{thm:main}.
In Section~\ref{sec:g2}, we follow the same procedures to give a~proof of Theorem~\ref{thm:main2}.

\section{Background}\label{sec:background}
\subsection*{Notation}
We will use the following notation throughout the paper:
\begin{itemize}\itemsep=0pt
\item $\g$: simple Lie algebra over $\C$ of rank $r$,
\item $\h$: Cartan subalgebra of $\g$,
\item $I=\{1,\dots, r\}$: index set for the Dynkin diagram of $\g$,
\item $\alpha_a$, $a\in I$: simple root,
\item $h_a$, $a\in I$: simple coroot,
\item $\omega_a$, $a\in I$: fundamental weight,
\item $C = (C_{ab})_{a,b\in I}$: Cartan matrix with $C_{ab}=\alpha_b(h_a)$,
\item $P=\bigoplus_{a\in I}\Z\omega_a$: weight lattice,
\item $P^{+}\subseteq P$: set of dominant integral weights,
\item $Q$: root lattice,
\item $\rho = \sum_{a\in I}\omega_a$: Weyl vector,
\item $\theta\in Q$: highest root,
\item $\posR$: set of positive roots,
\item $\h_{\R}^{*}:=\bigoplus_{a\in I}\R\omega_a$,
\item $(\cdot,\cdot)\colon \h_{\R}^{*}\times \h_{\R}^{*} \to \R$: $\R$-bilinear form induced from the Killing form with $(\theta,\theta)=2$,
\item $\Z[P]$: integral group ring of $P$ (which is the same as the ring \smash{$\Z\big[{\rm e}^{\pm \omega_j}\big]_{j\in I}$} of Laurent polynomials in ${\rm e}^{\omega_j}$),
\item $\K:=\C({\rm e}^{\omega_j})_{j\in I}$: field of rational functions in ${\rm e}^{\omega_j}$ with coefficients in $\C$,
\item $t_a: =(\theta,\theta)/(\alpha_a,\alpha_a)\in \{1,2,3\}$,
\item $[\alpha]_a\in \Z (\alpha \in Q,a\in I)$: coefficients in the sum \smash{$\alpha = \sum_{a\in I}[\alpha]_a\alpha_a$},
\item $s_a$: simple reflection acting on $\h_{\R}^{*}$ by $s_a(\la) = \la - \la(h_a)\alpha_a$,
\item $W$: Weyl group generated by $\{s_a \mid a\in I\}$,
\item $W_{\la}$, $\la\in P$: isotropy subgroup of $W$ fixing $\la$,
\item $W_{J}$, $J\subseteq I$: standard parabolic subgroup of $W$ generated by $\{s_a \mid a\in J\}$,
\item $L(\la)$, $\la \in P^{+}$: irreducible highest weight representation of $\uqg$,
\item $\chi(V)\in \Z[P]$: character of a~finite-dimensional $\uqg$-module $V=\bigoplus V_{\la}$ with weight spaces~$V_{\la}$, i.e., \smash{$\chi(V)=\sum_{\la\in P}(\dim V_{\la}){\rm e}^{\la}$}.
\end{itemize}

\subsection{Some properties of characters of KR modules}
KR modules form a~family of irreducible finite-dimensional representations of the quantum affine algebra $\uqghat$, where $q\in \C$ is not a~root of unity. For every \smash{$(a,m,u) \in I\times \Z_{\geq 0}\times \C^{\times}$}, there exists a~corresponding KR module \smash{$W^{(a)}_{m}(u)$}. By restriction, we obtain a~finite-dimensional $\uqg$-module \smash{$\res \Wam(u)$}, which can be denoted by \smash{$\res \Wam$} since its isomorphism class as a~$\uqg$-module does not depend on the spectral parameter $u$.

Let \smash{$\Qam:=\chi\bigl(\rwma\bigr)$}. The $Q$-system
\begin{equation}\label{eq:Qsys}
\bigl(\Qam\bigr)^2 = Q^{(a)}_{m+1}Q^{(a)}_{m-1} +
\prod_{b : C_{ab}< 0} \prod_{k=0}^{-C_{a b}-1}
Q^{(b)}_{\lfloor\frac{C_{b a}m - k}{C_{a b}}\rfloor}, \qquad a\in I,\quad m \geq 1,
\end{equation}
is a~difference equation satisfied by the characters of KR modules. Nakajima~\cite{MR1993360} and Hernandez~\cite{MR2254805} proved that the $q$-characters of KR modules satisfy the $T$-system, from which we obtain the $Q$-system~\eqref{eq:Qsys} by ignoring the spectral parameter.

In~\cite{MR3948941}, we studied a~linear recurrence relation with constant coefficients satisfied by the sequence \smash{$\bigl(\Qam\bigr)_{m=0}^{\infty}$}. We can summarize its main properties in terms of its generating function~\smash{$\Qa(t) : = \sum_{m=0}^{\infty}\Qam t^m$} as follows.

\begin{Theorem}[{\cite[Theorem 1.1]{MR3948941}}]\label{thm:linQ}
Let $\g$ be a~simple Lie algebra and $a\in I$.
Assume that $(\g,a)$ is not one of the following exceptional cases: node $4$ in type $E_7$, and nodes $3$, $4$, $5$, $8$ in type $E_8$.
There exist $W$-invariant finite subsets $\La_a$ and $\La'_a$ of $P$ with the following properties:
\begin{enumerate}[label=$(\roman*)$, ref=$(\roman*)$]\itemsep=0pt
\item If we set $D^{(a)}(z) : = \prod_{\la\in \La_a}\bigl(1-{\rm e}^{\la}z\bigr)\prod_{\la\in \La_a'}\bigl(1-{\rm e}^{\la}z^{t_a}\bigr)$, then
\begin{equation}\label{eqn:Qgen}
N^{(a)}(z): = \Qa(z)D^{(a)}(z)
\end{equation}
is a~polynomial in $z$ with coefficients in $\Z[P]$ and $\deg N^{(a)}<\deg D^{(a)}$.
\item \label{thm:lineq0} $t_a\La_{a}\cap \La'_{a}=\varnothing$, where $t_a\La_{a}=\{t_a \la\mid \la \in \La_a\}$.
\item $\omega_a\in \La_a$.
\end{enumerate}
\end{Theorem}
The same form of the linear recurrence is expected to hold for the excluded exceptional nodes as well, although those cases are not needed for the applications considered in this paper.
Let us fix $\La_a$ and $\La'_a$ as in the Appendix of~\cite{MR3948941}. When $\g$ is simply-laced, $\La_a$ is simply the set of weights of the fundamental representation $L(\omega_a)$.

Note that~\ref{thm:lineq0} shows that $1/D^{(a)}(t)$ has only simple poles.
Then, from the partial fraction decomposition of $\Qa(t)=N^{(a)}(t)/D^{(a)}(t)$, we deduce that for each $(\la,\zeta,l)\in P\times \C^{\times} \times \Z_{>0}$ there exists $C\bigl(\Qa,\la,\zeta,l\bigr)\in \K\bigl({\rm e}^{\omega_j/{t_a}}\bigr)_{j\in I}$ such that
\begin{equation}\label{eq:QandC0}
\Qam =\sum_{(\la,\zeta,l)} C\bigl(\Qa,\la,\zeta,l\bigr)\zeta^m {\rm e}^{m \la/l}, \qquad m\in \Z_{\geq 0},
\end{equation}
and the coefficient $C\bigl(\Qa,\la,\zeta,l\bigr)$ is zero unless either
\begin{itemize}\itemsep=0pt
\item $(\la,\zeta,l)=(\la,1,1)$ with $\la\in \La_a$; or,
\item $(\la,\zeta,l)=(\la,\zeta,t_a)$ with $\la\in \La_a'$ and $\zeta^{t_a}=1$.
\end{itemize}
Due to the $W$-symmetry of $\Qam$, we have
\begin{equation}\label{eqn:RWsym}
w\bigl(C\bigl(\Qa,\la,\zeta,l\bigr)\bigr)= C\bigl(\Qa,w(\la),\zeta,l\bigr),\qquad w\in W.
\end{equation}

Let \smash{$C^{(a)}_{\la}:=C\bigl(\Qa,\la,1,1\bigr)$} for $\la\in \La_a$ and \smash{$C^{(a)}_{\la,\zeta}:=C\bigl(\Qa,\la,\zeta,t_a\bigr)$} for $\la\in \La_a'$ and $\zeta^{t_a}=1$. We can rewrite~\eqref{eq:QandC0} as
\begin{equation}\label{eq:QandC}
\Qam =\sum_{\la\in \La_a} C^{(a)}_{\la}{\rm e}^{m \la}+\sum_{\la\in \La_a'}\sum_{\zeta\colon \zeta^{t_a}=1} C^{(a)}_{\la,\zeta}\zeta^m {\rm e}^{m \la/t_a}.
\end{equation}
For $C^{(a)}_{\omega_a}$, we have an explicit product formula.
\begin{Theorem}\label{thm:MYformula}
For each $a\in I$,
\begin{equation}\label{eqn:MYformula}
C^{(a)}_{\omega_a}=\frac{1}{\prod_{\alpha\in \posR}\bigl(1-{\rm e}^{-\alpha}\bigr)^{[\alpha]_a}}.
\end{equation}
Here, $[\alpha]_a\in \Z_{\geq 0}$ denotes the coefficient in the sum $\alpha =\sum_{a\in I}[\alpha]_a\alpha_a$.
\end{Theorem}
We call~\eqref{eqn:MYformula} the Mukhin--Young formula, which was originally conjectured in~\cite{MR3217701}. It was proved in~\cite{MR4283572} except for two nodes in type $E_8$; these remaining cases are covered in the more general setting of~\cite{arXiv250100724}. We note that, in general, coefficients other than $C\bigl(\Qa,\la,1,1\bigr)$ do not seem to admit an expression as compact as~\eqref{eqn:MYformula}.

\subsection{Fermionic formula}
The fermionic formula, proposed by Kirillov and Reshetikhin~\cite{MR906858}, concerns the decomposition of a~tensor product of Kirillov--Reshetikhin modules into irreducible $\uqg$-modules. Let $\smash{\bigl(\nu^{(a)}_m\bigr){}_{a\in I,m \geq 1}}$ be a~family of non-negative integers such that \smash{$\nu^{(a)}_m$} is zero for all but finitely many~${(a,m)}$. Consider
\[
W = \bigotimes_{(a,m)}\bigl(\res W^{(a)}_m\bigr)^{\otimes \nu^{(a)}_m},
\]
 a~tensor product of Kirillov--Reshetikhin modules, and its decomposition into irreducible $\uqg$-representations
\[
W \cong \bigoplus_{\la\in P^{+}}m(W,\la)L(\la),\qquad m(W,\la)\in \Z_{\geq 0},
\]
where $L(\la)$ denotes an irreducible $\uqg$-representation with highest weight $\la$. The fermionic formula provides an explicit combinatorial description of the multiplicity $m(W,\la)$ in terms of~\smash{$\bigl(\nu^{(a)}_m\bigr)_{a\in I,m \geq 1}$} and $\la$. Since this formula is somewhat complicated and not essentially used in this paper, we refer the reader, for example, to~\cite{MR1745263} for its precise statement.

\subsection{Polyhedral formula}\label{ss:polyform}
Recall that a~(rational) polyhedron $P\subset \R^d$ is the set of solutions of a~finite number, say $l$, of linear inequalities with integer coefficients
$
P=\{\mathbf{x} \in \R^d\mid \langle \mathbf{a}_i,\mathbf{x}\rangle \leq \beta_i \text{ for }i=1,\dots, l\}
$
for some~${\mathbf{a}_i\in \Z^d}$ and $\beta_i\in \Z$.
Here $\langle \cdot,\cdot \rangle$ denotes the standard inner product on $\R^d$.
We may include linear equalities because $\langle \mathbf{a}_i,\mathbf{x}\rangle \leq \beta_i$ and $\langle \mathbf{-a}_i,\mathbf{x}\rangle \leq -\beta_i$ imply $\langle \mathbf{a}_i,\mathbf{x}\rangle = \beta_i$.

Let $\theta=\sum_{a\in I}c_a\alpha_a, c_a\in \Z_{>0}$ be the highest root of $\g$. Fix $a\in I$ such that $c_a\leq 2$.
For example, this condition holds for all classical types.
It is shown in~\cite{MR1836791} (see also~\cite{MR2238884}) that there exist positive integers $(b_j)_{j\in J_{a}}$, and dominant integral weights $(\la_j)_{j\in J_{a}}$ for some finite set $J_{a}$ such that
\[
\rwma=\bigoplus_{\mathbf{x}\in \Fam} L(\la_{\mathbf{x}})
\]
where
\smash{$\Fam=\bigl\{(x_j)_{j\in J_{a}}\mid \sum_{j\in J_{a}}b_jx_j=m, x_j\in \Z_{\geq 0}\bigr\}$}, and \smash{$\la_{\mathbf{x}}=\sum_{j\in J_{a}}x_j\la_j$} for each \smash{$\mathbf{x}\in \Fam$}.
Hence, \smash{$\Fam$} is the integer lattice points of a~rational polyhedron.
When $c_a>2$, we still expect to have a~similar formula, but now with multiplicity, of the form
\begin{equation}\label{KRdecpolyform}
\rwma
\stackrel{?}{=}
\bigoplus_{\mathbf{x}\in \Fam} p(\mathbf{x};m)L(\la_{\mathbf{x}}).
\end{equation}
We further expect that $p(\mathbf{x};m)$ is $\Z_{\geq 0}$-valued for $m\in \Z_{\geq 0}$, \smash{$\mathbf{x}\in \Fam$}, and can be written as step-polynomials defined piecewise on polyhedral regions.
A~\textit{step-polynomial} is a~$\QQ$-linear combination of products of functions of the form $\lfloor \langle \mathbf{a},\cdot \rangle+b \rfloor$ for some $\mathbf{a}\in \QQ^{d}$, $b\in \QQ$.
A polyhedral formula for the decomposition of KR modules will mean a~formula of the form~\eqref{KRdecpolyform}.

For example, consider $p(j_1,j_2,j_3,j_4;m)$ from Theorem~\ref{thm:main}.
Then $p = g_i$ on $Q_i$ ($i=1,2$), where
\begin{gather*}
Q_1 = \{(j_1,\dots, j_4,m)\in \Z_{\geq 0}^5 \mid j_1 + 2 j_2 + 2j_3 + j_4 \le m\},\qquad
g_1 = (1+j_3)(j_4+1),\\
Q_2 = \{(j_1,\dots, j_4,m)\in \Z_{\geq 0}^5 \mid j_1 + 2 j_2 + j_3 + j_4 \le m \le j_1 + 2 j_2 + 2j_3 + j_4\},\\
g_2 = (1+m-j_1-2j_2-j_3-j_4)(j_4+1).
\end{gather*}
In the $a=2$ case of type $G_2$, the multiplicity function in~\eqref{eqn:step-poly} demonstrates the necessity of the notion of step-polynomials.

As mentioned in the introduction, computations based on the fermionic formula can provide useful data for guessing the multiplicity function $p(\mathbf{x};m)$ in~\eqref{KRdecpolyform}. However, a~candidate formula for $p(\mathbf{x};m)$ still has to be guessed in advance; the method developed below is designed to verify such a~proposed formula for all $m\in \Z_{\geq 0}$. In the next section, we explain an approach for a~proof of~\eqref{KRdecpolyform}.

\section{Framework for proving polyhedral formula}\label{sec:framework}
Let \smash{$\Pam$} denote the character of the right-hand side of~\eqref{KRdecpolyform}, i.e.,
\[
\Pam=\sum_{\mathbf{x}\in \Fam}p(\mathbf{x};m)\chi(L(\la_{\mathbf{x}})).
\]
The letter $P$ is chosen from the word `polyhedron'.
When there is such a~polyhedral formula, consider its generating function
\smash{$
\Pa(t) : = \sum_{m=0}^{\infty}\Pam t^m$}.
Once we combine the discussion from Section~\ref{ss:polyform} and the Weyl character formula, we can obtain $\Pa(t)$ as a~rational function in $t$.

Then, \eqref{KRdecpolyform} is equivalent to the following identity between two rational functions
\begin{equation}\label{eq:QeqP}
\Pa(t) = \Qa(t).
\end{equation}

We can state the steps necessary to prove~\eqref{eq:QeqP} as follows:
\begin{enumerate}[label=(\roman*)]\itemsep=0pt
\item\label{algo:resloc} Prove that \smash{$\Pa$} has at most simple poles, and the set of poles of \smash{$\Pa$} is a~subset of the set of poles of \smash{$\Qa$}. Also make sure that the degree of the denominator of \smash{$\Pa$} is greater than that of its numerator.
\item\label{algo:rescompe} Prove the equality of the coefficients
\begin{equation}\label{eqn:CPCQ}
C\bigl(\Pa,\la,\zeta,l\bigr) = C\bigl(\Qa,\la,\zeta,l\bigr)
\end{equation}
when $(\la,\zeta,l)$ belongs to one of the following cases:
\begin{itemize}\itemsep=0pt
\item $(\la,\zeta,l)=(\la,1,1)$ with $\la\in \La_a$;
\item $(\la,\zeta,l)=(\la,\zeta,t_a)$ with $\la\in \La_a'$ and $\zeta^{t_a}=1$.
\end{itemize}
\end{enumerate}

As \smash{$\Pa$} is a~rational function in $t$, we can effectively attack step~\ref{algo:resloc}, for example, using the partial fraction decomposition method as we do in Section~\ref{subsec:resloc}.

Recall that we already know explicitly where the poles of \smash{$\Qa$} are located and that these poles are always simple. After settling step~\ref{algo:resloc}, we can deduce that \smash{$\Pam$} can be written in the form
\[
\Pam =\sum_{(\la,\zeta,l)} C\bigl(\Pa,\la,\zeta,l\bigr)\zeta^m {\rm e}^{m \la/l}
\]
with \smash{$C\bigl(\Pa\!,\la,\zeta,l\bigr)\!\!\in\! \K\bigl({\rm e}^{\omega_j/{t_a}}\bigr)_{\!j\in I}$}, which vanishes unless the non-vanishing conditions for $C\bigl(\smash{\Qa}\!,\allowbreak\la,\zeta,l\bigr)$, stated after~\eqref{eq:QandC0}, are satisfied.

Step~\ref{algo:rescompe} is equivalent to showing that \smash{$\Pa$} and \smash{$\Qa$} have the same residues at their poles, which are simple at most.
Since \smash{$\Pa$} is $W$-invariant as each \smash{$\Pam$} is $W$-invariant, the coefficients~${C\bigl(\Pa,\la,\zeta,l\bigr)}$ follow the same $W$-symmetry as $C\bigl(\Qa,\la,\zeta,l\bigr)$ in~\eqref{eqn:RWsym}.
Hence, it is enough to consider weights in $P^{+}$ because every element of $P$ has a~unique element in $P^{+}$ in its $W$-orbit.

Once \smash{$\Pam$} is explicitly given as in~\eqref{KRdecpolyform}, it is more or less straightforward to compute \smash{$\Pa$} and $C\bigl(\Pa,\la,\zeta,l\bigr)$.
For computing $C\bigl(\Qa,\la,\zeta,l\bigr)$, we can use the $Q$-system~\eqref{eq:Qsys} to derive a~system of equations among such coefficients, and find its solution, usually in a~recursive way from known ones.

Suppose that we already have proved~\ref{algo:resloc}.
Then another way to finish the proof of $\smash{\Pa(t)} = \smash{\Qa(t)}$ is by showing \smash{$\Qam = \Pam$} for $m=0,1,\dots, |\La_a|+t_a |\La_a'|-1$ since now we know that both sequences satisfy the same linear recurrence relation of order $|\La_a|+t_a |\La_a'|$.
Although it is a~purely mechanical task to check \smash{$\Qam = \Pam$} for given $m$ using the fermionic formula, we have found that it is still computationally challenging for large $m$.

By considering \smash{$C\bigl(\Pa,\la,\zeta,l\bigr)$} and \smash{$C\bigl(\Qa,\la,\zeta,l\bigr)$}, we localize the problem in the sense that we are looking at a~single pole at a~time, and thus obtain further simplifications.
Another advantage is that we are able to break the original problem into smaller and more tractable pieces.

In a~nutshell, it is possible to check both~\ref{algo:resloc} and~\ref{algo:rescompe} algorithmically.
In the following sections, we follow this strategy to prove polyhedral formulas in type $F_4$ and in type $G_2$ where we discuss some practical issues in our method in detail.

\begin{Remark}\label{rem:szlemma}
In carrying out step~\ref{algo:rescompe},
we are eventually led to the problem of showing that a~sum of many multivariate rational functions actually vanishes.
After removing the common denominator, this is essentially a~problem about multivariate polynomials.

Polynomial identity testing (PIT) is one of the important problems
at the intersection of algebra and computer science.
It might be useful to add a~few words about a~well-known result in the area.
Let $F$ be a~field, $P\in F[x_1,x_2,\ldots,x_n]$ a~non-zero polynomial, $S$ a~finite non-empty subset of $F$.
The Schwartz--Zippel lemma gives an explicit upper bound on the size of~${\{(r_1,\dots,r_n)\in S^{n} \mid P(r_1,\ldots,r_n)=0\}}$
in terms of $|S|$ and the degree of $P$.
We refer the reader to~\cite{MR3788163} and the references therein for a~precise statement.
\end{Remark}

\begin{Remark}
We know that \smash{$C\bigl(\Qa,\la,\zeta,l\bigr)$} is invariant under $W_{\la}$ from~\eqref{eqn:RWsym}. When we explicitly compute \smash{$C\bigl(\Pa,\la,\zeta,l\bigr)$}, it is given as a~sum over $W_{\la}$; see~\eqref{eqn:D2la} for an example. Thus we can regard~\eqref{eqn:CPCQ} as a~summation formula over $W_{\la}$ for \smash{$C\bigl(\Qa,\la,\zeta,l\bigr)$}.
\end{Remark}

The accompanying \textsc{Mathematica} notebook file for some computer calculations in this paper is available at \url{https://github.com/chlee-0/KR-polyhedral-formula}.

\section[Proof of Theorem 1.1 for type F\_4]{Proof of Theorem~\ref{thm:main} for type $\boldsymbol{F_4}$}\label{sec:f4}
Let $\g$ be a~simple Lie algebra of type $F_4$. When $a=1$ or $4$, there is a~known polyhedral formula for \smash{$\res \Wam$}. The formula for $a=1$ is given in~\eqref{Qonem} and will be used later. The main goal of this section is to prove Theorem~\ref{thm:main}, namely, the polyhedral formula for $a=2$.

For $m\in \Z_{\geq 0}$, let \smash{$\PP{2}{m}$} be the character of the right-hand side of~\eqref{eqn:sumef} and $\smash{\mathcal{P}^{(2)}(t)} = \smash{\sum_{m=0}^{\infty}\PP2m t^m}$. By following the strategy outlined in Section~\ref{sec:framework}, we will show that
\begin{equation}\label{eqn:QeqP}
\Q2m =\PP{2}{m},
\end{equation}
that is,
$
\mathcal{Q}^{(2)}(t)=\mathcal{P}^{(2)}(t)
$
as rational functions in $t$.

Before turning to proofs, we present Table~\ref{tab:F4-lambda-sets} for $\La_a$ and $\La_a'$ from~\cite{MR3948941}. Since they are $W$-invariant, they are given as a~disjoint union of $W$-orbits of elements of $P^+$.
\begin{table}[ht]\centering
$
\begin{array}{c|c|c}
 a~& \La_a\cap P^+ & \La'_a\cap P^+ \\
\hline
 1 & 0,\,\omega_1 & \cdot \\
 2 & 0,\,\omega_1,\,\omega_2,\,2 \omega_4 & \cdot \\
 3 & \omega_3 & 0,\omega_1,\,2 \omega_1,\,\omega_2,\,2 \omega_4,\,\omega_1+2 \omega_4 \\
 4 & \omega_4 & 0,\,\omega_1 \\
\end{array}
$
\caption{The sets $\La_a$ and $\La'_a$ in type $F_4$.}
\label{tab:F4-lambda-sets}
\end{table}
Throughout the section, \smash{$\Delta = \sum_{w\in W} (-1)^{\ell(w)}{\rm e}^{w(\rho)}$} denotes the Weyl denominator and $(\la_1,\la_2,\la_3,\la_4)=(2\omega_4,2\omega_3,\omega_2,\omega_1)$.
We often need to explicitly deal with $W$ or its subgroups.
One may refer to~\cite{MR1073413} for an algorithm to find the Weyl orbit of a~weight or a~minimal coset representative for a~coset of a~standard parabolic subgroup of a~Weyl group.

\subsection{Step~\ref{algo:resloc}}\label{subsec:resloc}
From~\eqref{eqn:Qgen}, it is clear that \smash{$\mathcal{Q}^{(2)}(t)$} can have only poles of order at most 1 and they can only be found at $t={\rm e}^{-\la}$, $\la\in \La_2$. Let us find the poles of \smash{$\mathcal{P}^{(2)}(t)$}. To write \smash{$\mathcal{P}^{(2)}(t)$} explicitly, define a~sequence $\{a_m\}_{m=0}^{\infty}$ by
\[
a_m=\sum_{
	\substack{
		j_1 + 2 j_2 + j_3 + j_4 \le m\\
		j_1,j_2,j_3,j_4 \in \Z_{\geq 0}
	}
	}
p(j_1,j_2,j_3,j_4;m) x_1^{j_1}x_2^{j_2}x_3^{j_3}x_4^{j_4},
\]
whose generating function is
\[
\sum_{m=0}^{\infty}a_mt^m=\frac{1}{(1-t) (1-t x_1) \bigl(1-t^2 x_2\bigr) (1-t x_3) \bigl(1-t^2 x_3\bigr) (1-t x_4)^2}.
\]
We note that this sum can be found using only elementary facts about geometric series.

By combining this with the Weyl character formula, \smash{$\mathcal{P}^{(2)}(t)$} can be written as
\begin{equation}\label{eqn:Psum}
\frac{1}{\Delta}\sum_{w\in W} \frac{(-1)^{\ell(w)}{\rm e}^{w(\rho)}}{(1-t)\bigl(1-{\rm e}^{w(\la _1)}t\bigr)\bigl(1-{\rm e}^{w(\la _2)}t^2\bigr)\bigl(1-{\rm e}^{w(\la_3)}t\bigr)\bigl(1-{\rm e}^{w(\la _3)}t^2\bigr)\bigl(1-{\rm e}^{w(\la _4)}t\bigr)^2}.
\end{equation}
At this point, it is not entirely clear whether \smash{$\mathcal{P}^{(2)}(t)$} has only simple poles at $t={\rm e}^{-\la}$, $\la\in \La_2$ or not. For example, $\mathcal{P}^{(2)}(t)$ may have a~double pole at \smash{$t={\rm e}^{-w(\la_4)}$}.

Consider the partial fraction decomposition of a~summand in~\eqref{eqn:Psum}
\begin{gather*}
\frac{1}{(1-t)\bigl(1-{\rm e}^{w(\la _1)}t\bigr)\bigl(1-{\rm e}^{w(\la _2)}t^2\bigr)\bigl(1-{\rm e}^{w(\la_3)}t\bigr)\bigl(1-{\rm e}^{w(\la _3)}t^2\bigr)\bigl(1-{\rm e}^{w(\la _4)}t\bigr)^2} \\
\qquad=
\frac{w\bigl(d^{(2)}_{0}\bigr)}{1-t}
+\frac{w\bigl(d^{(2)}_{\la_1}\bigr)}{1-{\rm e}^{w(\la_1)}t}
+\frac{w\bigl(d^{(2)}_{\la_3}\bigr)}{1-{\rm e}^{w(\la_3)}t}
+\frac{w\bigl(d^{(2)}_{\la_4}\bigr)}{1-{\rm e}^{w(\la_4)}t} \\
\phantom{\qquad=}{}
+\frac{w\bigl(e^{(2)}_{\la_2}\bigr)}{1-{\rm e}^{w(\la_2)}t^2}
+\frac{w\bigl(e^{(2)}_{\la_3}\bigr)}{1-{\rm e}^{w(\la_3)}t^2}
+\frac{w\bigl(e^{(2)}_{\la_4}\bigr)}{\bigl(1-{\rm e}^{w(\la_4)}t\bigr)^2},
\end{gather*}
where \smash{$d^{(2)}_{\la}, e^{(2)}_{\la_4}\in \K$}, and \smash{$e^{(2)}_{\la_2},e^{(2)}_{\la_3}\in \K[t]$} are polynomials of degree at most~1, uniquely determined by this form of decomposition in the case $w=1$.
For example,
\begin{equation}\label{Dtwozero}
d^{(2)}_{0}=
\frac{1}{\bigl(1-{\rm e}^{\la _1}\bigr) \bigl(1-{\rm e}^{\la _2}\bigr) \bigl(1-{\rm e}^{\la _3}\bigr)^2 \bigl(1-{\rm e}^{\la _4}\bigr)^2}.
\end{equation}
Because these expressions are long but easy to find, we do not write them here; one can refer to the accompanying file for an explicit description.
\begin{Proposition}\label{prop:pfd}
We have
\[
\mathcal{P}^{(2)}(t) =
\frac{1}{\Delta}\sum_{w\in W}(-1)^{\ell(w)}{\rm e}^{w(\rho)}
\left(
\frac{w\bigl(d^{(2)}_{0}\bigr)}{1-t}
+\frac{w\bigl(d^{(2)}_{\la_1}\bigr)}{1-{\rm e}^{w(\la _1)}t}
+\frac{w\bigl(d^{(2)}_{\la_3}\bigr)}{1-{\rm e}^{w(\la _3)}t}
+\frac{w\bigl(d^{(2)}_{\la_4}\bigr)}{1-{\rm e}^{w(\la _4)}t}
\right).
\]
\end{Proposition}

\begin{proof}
It is sufficient to show that for $\la\in \{\la_2,\la_3,\la_4\}$,
\begin{equation}\label{sumEvanishing}
\sum_{w\in W_{\la}}(-1)^{\ell(w)}{\rm e}^{w(\rho)}w\bigl(e^{(2)}_{\la}\bigr)=0.
\end{equation}
We may use computers to verify this directly. Below, we will explain how to reduce the amount of calculation to check~\eqref{sumEvanishing}. While this reduction is not essential as long as we focus on type~$F_4$ whose Weyl group is manageable in size, it might be useful for treating a~similar vanishing sum over bigger groups in other types.

For $\la=\la_2=2\omega_3$, $W_{\la}=W_{\{1,2,4\}}$. The parabolic subgroup $W_{\{1,2\}}$ of $W_{\{1,2,4\}}$ satisfies
\[
\sum_{w\in W_{\{1,2\}}}(-1)^{\ell(w)}{\rm e}^{w(\rho)}w\bigl(e^{(2)}_{\la_2}\bigr)=0,
\]
which implies the vanishing of the sum over $W_{\{1,2,4\}}$ since
\eqref{sumEvanishing} can be written as
\[
\sum_{w'\in W'}(-1)^{\ell(w')}w'\left(\sum_{w\in W_{\{1,2\}}}(-1)^{\ell(w)}{\rm e}^{w(\rho)}w\bigl(e^{(2)}_{\la_2}\bigr)\right),
\]
where $W'$ is the set of minimal coset representatives for cosets in $W_{\{1,2,4\}}/W_{\{1,2\}}$.

Similarly, when $\la=\la_3=\omega_2$, we have $W_{\la}=W_{\{1,3,4\}}$ and the following sum over the parabolic subgroup $W_{\{1,3\}}$ vanishes:
\[
\sum_{w\in W_{\{1,3\}}}(-1)^{\ell(w)}{\rm e}^{w(\rho)}w\bigl(e^{(2)}_{\la_3}\bigr)=0.
\]

When $\la=\la_4=\omega_1$, we have not found any proper parabolic subgroup of $W_{\la}=W_{\{2,3,4\}}$, over which the sum vanishes. However, if we let
\[
E^{(2)}_{\la_4}:=\sum_{w\in W_{\{2,4\}}}(-1)^{\ell(w)}{\rm e}^{w(\rho)}w\bigl(e^{(2)}_{\la_4}\bigr),
\]
then the left-hand side of~\eqref{sumEvanishing} becomes
\[
\sum_{w\in W''}(-1)^{\ell(w)} w\bigl(E^{(2)}_{\la_4}\bigr),
\]
where $W''$ is the set of minimal coset representatives for cosets in $W_{\{2,3,4\}}/W_{\{2,4\}}$. The size of~$W''$ is 12, and it is possible to partition this set into 6 pairs of distinct elements so that the contribution from each pair to the above sum is zero. In other words, for each $w_1\in W''$, there exists $w_2\in W''$, $w_1\neq w_2$ such that
\[
(-1)^{\ell(w_1)} w_1\bigl(E^{(2)}_{\la_4}\bigr)+(-1)^{\ell(w_2)}w_2\bigl(E^{(2)}_{\la_4}\bigr)=0.\tag*{\qed}
\] \renewcommand{\qed}{}
\end{proof}

This proposition immediately implies the following.
\begin{Proposition}\label{prop:Ppoles}
The rational function $\mathcal{P}^{(2)}(t)$ has only simple poles, possibly at $t={\rm e}^{-\la}$, $\la\in \La_2$ and no other poles.
\end{Proposition}

Now we know that the poles of $\mathcal{Q}^{(2)}(t)$ and $\mathcal{P}^{(2)}(t)$ can only appear at $t= {\rm e}^{-\la}$, $\la \in \La_2$.

\subsection[Step (ii)]{Step~\ref{algo:rescompe}}\label{ss:rescomp}
It remains to carry out the second step of our strategy from Section~\ref{sec:framework} to prove~\eqref{eqn:QeqP}.
As shown in Table~\ref{tab:F4-lambda-sets}, $\La'_2$ is empty. Hence we have to show
$
C\bigl(\mathcal{P}^{(2)},\la,1,1\bigr) = C\bigl(\mathcal{Q}^{(2)},\la,1,1\bigr)
$
for $\la \in \La_2\cap P^+=\{0,\omega_1,\omega_2,2 \omega_4\}$. We first explain how to compute both sides, and then check their equality. Let us write \smash{$C^{(2)}_{\la}=C\bigl(\mathcal{Q}^{(2)},\la,1,1\bigr)$} and \smash{$D^{(2)}_{\la}= C\bigl(\mathcal{P}^{(2)},\la,1,1\bigr)$}.

\subsubsection*{How to compute $\boldsymbol{C^{(2)}_{\la}}$}
Let us explain how to compute \smash{$C^{(2)}_{\la}$}.
Recall the $Q$-system relation~\eqref{eq:Qsys}
$
\smash{Q^{(2)}_{m} = \bigl(Q^{(1)}_{m}\bigr)^2 }- \smash{Q^{(1)}_{m-1}Q^{(1)}_{m+1}}
$
for $a=1$.
By rewriting this relation using~\eqref{eq:QandC}, we obtain an expression for~\smash{$C^{(2)}_{\la}$} in terms of~\smash{$C^{(1)}_{\mu}$}, $\mu \in \La_1$
\begin{equation}\label{eqn:C2C1}
C^{(2)}_{\la} = \sum_{(\nu,\mu) \in S_{\la}} C^{(1)}_{\nu}C^{(1)}_{\mu}(1-{\rm e}^{\nu-\mu}),
\end{equation}
where
\smash{$S_{\la} : =\{(\nu,\mu)\in \La_1\times \La_1\mid \nu \neq \mu,\, \nu+\mu = \la\}$}. To handle~\eqref{eqn:C2C1} explicitly, we need a~way~to compute \smash{$C^{(1)}_{0}$} and \smash{$C^{(1)}_{\omega_1}$}. And these are all we need to
find \smash{$C^{(1)}_{\mu}$}, \smash{$\mu \in \La_1$} because \smash{$C^{(1)}_{\la}$} with non-zero $\la\in W\omega_1$ is given by \smash{$C^{(1)}_{w(\omega_1)} = w\bigl(C^{(1)}_{\omega_1}\bigr)$}.
\begin{Lemma}\label{lem:C10}
We have
\begin{equation}\label{C10}
C^{(1)}_{0} = \frac{\sum_{w\in W}(-1)^{\ell(w)}{\rm e}^{w(\rho)}/\bigl(1-{\rm e}^{w(\omega_1)}\bigr)}{\Delta},
\end{equation}
and
\[
C^{(1)}_{\omega_1} = 1/\prod_{\alpha\in \posR}\bigl(1-{\rm e}^{-\alpha}\bigr)^{[\alpha]_1}.
\]
\end{Lemma}

\begin{proof}
Note that $C^{(1)}_{\omega_1}$ is given by Theorem~\ref{thm:MYformula}, the Mukhin--Young formula. To find $C^{(1)}_{0}$, we can exploit the known polyhedral formula from~\cite{MR1836791}
\begin{equation}\label{Qonem}
Q^{(1)}_m = \sum_{k=0}^{m}\chi(L(k\omega_1)).
\end{equation}
By the Weyl character formula,
\begin{align*}
\mathcal{Q}^{(1)}(t) & = \frac{1}{\Delta}\sum_{w\in W} \frac{(-1)^{\ell(w)}{\rm e}^{w(\rho)}}{(1-t)\bigl(1-t {\rm e}^{w(\omega_1)}\bigr)} = \frac{1}{\Delta}\sum_{w\in W}(-1)^{\ell(w)}{\rm e}^{w(\rho)}\left(\frac{w\bigl(d^{(1)}_{0}\bigr)}{1-t}+\frac{w\bigl(d^{(1)}_{\omega_1}\bigr)}{1-t {\rm e}^{w(\omega_1)}}\right),
\end{align*}
where
\smash{$
d^{(1)}_{0} = \frac{1}{1-{\rm e}^{\omega_1}}$} and \smash{$
 d^{(1)}_{\omega_1} = \frac{1}{1-{\rm e}^{-\omega_1}}$}.
Therefore,
\[
C^{(1)}_{0} = \frac{\sum_{w\in W}(-1)^{\ell(w)}{\rm e}^{w(\rho)}w\bigl(d^{(1)}_{0}\bigr)}{\Delta}.\tag*{\qed}
\] \renewcommand{\qed}{}
\end{proof}

\subsubsection*{How to compute $\boldsymbol{D^{(2)}_{\la}}$}
Recall that we have
\smash{$
\PP2m =\sum_{\la\in \La_2} D^{(2)}_{\la}{\rm e}^{m \la}$}.
We only need to consider $\la \in \La_2\cap P^+$.
By~Proposi\-tion~\ref{prop:pfd}, we can write \smash{$D^{(2)}_{\la}$} as
\begin{equation}\label{eqn:D2la}
D^{(2)}_{\la} = \frac{\sum_{w\in W_{\la}}(-1)^{\ell(w)}{\rm e}^{w(\rho)}w(d^{(2)}_{\la})}{\Delta}.
\end{equation}
For a~dominant weight $\la$, $W_{\la}$ is a~standard parabolic subgroup of $W$. Once we enumerate the elements of $W_{\la}$, it is straightforward to compute \smash{$D^{(2)}_{\la}$}. Of course, it becomes computationally easier to manipulate~\eqref{eqn:D2la} when $W_{\la}$ is a~proper subgroup of $W$. In this sense, the most difficult case arises when $\la=0$.

Now we can compute both~\eqref{eqn:C2C1} and~\eqref{eqn:D2la} and thus, are ready to check \smash{$C^{(2)}_{\la} = D^{(2)}_{\la}$} for $\la \in \{0,\omega_1,\omega_2,2 \omega_4\}$. To use~\eqref{eqn:C2C1}, we need $S_{\la} =\{(\nu,\mu)\in \La_1\times \La_1\mid \nu \neq \mu, \nu+\mu = \la\}$, as described below. The cases of $\la \in \{\omega_1,\omega_2,2 \omega_4\}$ do not bring much difficulty, and a~computer can easily simplify \smash{$C^{(2)}_{\la} - D^{(2)}_{\la}$} and return zero.
We give further comments on the $\la=0$ case, which is the most difficult one.

\subsubsection[lambda=omega\_1 case]{$\boldsymbol{\la=\omega_1}$ case}
In this case,
\begin{align*}
S_{\omega_1}={}&
\{
 (0,\omega _1), (\omega _1,0),
 (2 \omega _1-\omega _2,\omega _2-\omega _1), (\omega _2-\omega _1,2 \omega _1-\omega _2),
 (\omega _1+\omega _2-2 \omega _3,2 \omega _3-\omega _2), \\
& (2 \omega _3-\omega _2,\omega _1+\omega _2-2 \omega _3),
 (\omega _2-2 \omega _4,\omega _1-\omega _2+2 \omega _4), (\omega _1-\omega _2+2 \omega _4,\omega _2-2 \omega _4), \\
& (\omega _1-\omega _2+2 \omega _3-2 \omega _4,\omega _2-2 \omega _3+2 \omega _4), (\omega _2-2 \omega _3+2 \omega _4,\omega_1-\omega _2+2 \omega _3-2 \omega _4)
\}.
\end{align*}

\subsubsection[lambda=omega\_2 case]{$\boldsymbol{\la=\omega_2}$ case}
In this case, $S_{\omega_2} = \{(\omega _1,\omega _2-\omega _1), (\omega _2-\omega _1,\omega _1)\}$.
Thus~\eqref{eqn:C2C1} gives
\[
C^{(2)}_{\omega_2} = C_{\omega _1}^{(1)} C_{\omega _2-\omega _1}^{(1)}\bigl(2-{\rm e}^{2\omega_1-\omega_2}-{\rm e}^{\omega_2-2\omega_1}\bigr).
\]
In fact, this identity is a~special case of \cite[Corollary~3.5]{MR4283572}.

\subsubsection[lambda=2omega\_4 case]{$\boldsymbol{\la=2\omega_4}$ case}
In this case,
\begin{align*}
S_{2\omega_4} ={}&
\{
 (\omega _1,2 \omega _4-\omega _1), (2 \omega _4-\omega _1,\omega _1),
 (\omega _2-\omega _1,\omega _1-\omega _2+2 \omega _4),
 (\omega _1-\omega _2+2 \omega _4,\omega _2-\omega _1),\\
 &(2 \omega _3-\omega _2,\omega _2-2 \omega _3+2 \omega _4), (\omega _2-2 \omega _3+2 \omega _4,2 \omega _3-\omega _2)
\}.
\end{align*}

\subsubsection[lambda=0 case]{$\boldsymbol{\la=0}$ case}
Our goal is to check whether \smash{$C^{(2)}_{0}-D^{(2)}_{0}$} is actually zero, but this calculation is not quite straightforward as before, since it is a~sum of many rational functions.
\eqref{eqn:C2C1} can be rewritten as
\[
C^{(2)}_{0} = \sum_{\la\in W\omega_1}C^{(1)}_{\la}C^{(1)}_{-\la}\bigl(1-{\rm e}^{2\la}\bigr),
\]
and $|W\omega_1| = 24$. On the other hand, \smash{$D^{(2)}_{0}$} involves an alternating sum of orbits of~\eqref{Dtwozero} over the entire Weyl group $W$ and hence, it is obtained by adding $|W|=1152$ rational functions in~${{\rm e}^{\omega_1},\dots,{\rm e}^{\omega_4}}$.

It is slightly better to work with \smash{$\Delta C^{(2)}_{0}$} and \smash{$\Delta D^{(2)}_{0}$} to simplify their denominators. Let us consider
\[
\Delta D^{(2)}_{0}=\sum_{w\in W}(-1)^{\ell(w)}{\rm e}^{w(\rho)}w\bigl(d^{(2)}_{0}\bigr).
\]
We can rewrite the above as
\[
\Delta D^{(2)}_{0}=
\sum_{w'\in W^{\{1,3,4\}}}
(-1)^{\ell(w')}
w'(f_{w'}),
\]
where $W^{\{1,3,4\}}$ denotes the set of minimal coset representatives of cosets in $W/W_{\{1,3,4\}}$,
and
\[f_{w'}=
(-1)^{\ell(w')}
w'\left(\sum_{w\in W_{\{1,3,4\}}}(-1)^{\ell(w)}{\rm e}^{w(\rho)}w\bigl(d^{(2)}_{0}\bigr)\right),\qquad w'\in W^{\{1,3,4\}}
.\]
Note that the size of \smash{$W^{\{1,3,4\}}$} is 96.

In actual calculations, we treat the above expressions as rational functions of variables ${x_i={\rm e}^{\alpha_i}}$ with coefficients in $\QQ$.
Hence \smash{$\Delta C^{(2)}_{0}$} and \smash{$-\Delta D^{(2)}_{0}$} produce $24+96=120$ rational functions in~$x_i$, and we can simply try to sum them up to get zero.
On a~standard desktop computer, it took about 500 seconds to complete this calculation.
When we computed with 120 polynomials, after multiplying a~suitable common factor to these rational functions, it took about 1400 seconds to get zero.

\section[Proof of Theorem~1.2 for type G\_2]{Proof of Theorem~\ref{thm:main2} for type $\boldsymbol{G_2}$}\label{sec:g2}
Let $\g$ be of type $G_2$.
We record the sets $\Lambda_a$ and $\Lambda'_a$ as given in \cite[Appendix B]{MR3948941}
\[
\begin{array}{c|c|c} a~& \Lambda_a\cap P^+ & \Lambda'_a\cap P^+\\
\hline
1 & 0,\,\omega_1 & \cdot\\
2 & \omega_2 & 0,\,\omega_1
\end{array}
\]

Throughout this section, we assume that the polyhedral formula for \smash{$\res W^{(a)}_m$} is unknown in both cases $a=1$ and $a=2$, in order to illustrate our method.
The argument in this section gives a~proof for both cases simultaneously.
Let us first consider the case $a=1$.
Set~\smash{$\mathcal{P}^{(1)}(t)=\sum_{m=0}^{\infty}P^{(1)}_m t^m$}. Then
since the decomposition of \smash{$\res W^{(1)}_m$} agrees with~\eqref{Qonem},
we can reuse the computation in Lemma~\ref{lem:C10} to write \smash{$\mathcal{P}^{(1)}(t)$} as
\begin{align}
\mathcal{P}^{(1)}(t)&=
\frac{1}{\Delta}\sum_{w\in W} \frac{(-1)^{\ell(w)}{\rm e}^{w(\rho)}}{(1-t)\bigl(1-{\rm e}^{w(\omega_1)}t\bigr)}\nonumber\\
&=\frac{1}{\Delta}\sum_{w\in W}(-1)^{\ell(w)}{\rm e}^{w(\rho)}\left(\frac{w\bigl(d^{(1)}_{0}\bigr)}{1-t}+\frac{w\bigl(d^{(1)}_{\omega_1}\bigr)}{1-t {\rm e}^{w(\omega_1)}}\right).\label{eqn:PsumG2a1}
\end{align}
By comparing coefficients of $t^m$ in~\eqref{eqn:PsumG2a1}, we can rewrite \smash{$P^{(1)}_m$} in the form
\begin{equation}\label{eqn:pp1}
P^{(1)}_m=\sum_{\lambda\in \Lambda_1}D^{(1)}_{\lambda}{\rm e}^{m\lambda},
\qquad m\in \Z_{\geq 0},
\end{equation}
where
\smash{$
D^{(1)}_{0}:=
\frac{1}{\Delta}\sum_{w\in W}
(-1)^{\ell(w)}{\rm e}^{w(\rho)}w\bigl(d^{(1)}_{0}\bigr)$},
and
\[
D^{(1)}_{\lambda}:=
\frac{1}{\Delta}\sum_{\substack{w\in W\\ w(\omega_1)=\lambda}}
(-1)^{\ell(w)}{\rm e}^{w(\rho)}w\bigl(d^{(1)}_{\omega_1}\bigr),
\qquad \lambda\in W\omega_1.
\]
with \smash{$d^{(1)}_{0}=\frac{1}{1-{\rm e}^{\omega_1}}$} and \smash{$d^{(1)}_{\omega_1}=\frac{1}{1-{\rm e}^{-\omega_1}}$}.

For $m\in \Z_{\geq 0}$, let \smash{$\PP{2}{m}$} be the character of the right-hand side of~\eqref{eqn:sumeg} and set $\smash{\mathcal{P}^{(2)}(t) }= \smash{\sum_{m=0}^{\infty}\PP2m t^m}$.
We now write \smash{$\mathcal{P}^{(2)}(t)$} explicitly. To this end, define a~sequence $\{a_m\}_{m=0}^{\infty}$ by
\[
a_m=\sum_{
	\substack{
		j_0 + 3j_1 + j_2 = m\\
		j_0,j_1,j_2 \in \Z_{\geq 0}
	}
	}
p(j_0,j_1,j_2;m) x_1^{j_1}x_2^{j_2},
\]
whose generating function is
\[
\sum_{m=0}^{\infty}a_mt^m=\frac{1}{\bigl(1-t^3\bigr)\bigl(1-t^3 x_1\bigr)^2(1-t x_2)\bigl(1-t^2 x_2\bigr)}.
\]

By combining this with the Weyl character formula, \smash{$\mathcal{P}^{(2)}(t)$} can be written as
\begin{equation}\label{eqn:PsumG2}
\mathcal{P}^{(2)}(t)=
\frac{1}{\Delta \bigl(1-t^3\bigr)}\sum_{w\in W} \frac{(-1)^{\ell(w)}{\rm e}^{w(\rho)}}{\bigl(1-{\rm e}^{w(\omega_1)}t^3\bigr)^2\bigl(1-{\rm e}^{w(\omega_2)}t\bigr)\bigl(1-{\rm e}^{w(\omega_2)}t^2\bigr)}.
\end{equation}

Then for $w\in W$ we have
\begin{gather}
\frac{1}{\bigl(1-t^3\bigr)\bigl(1-{\rm e}^{w(\omega_1)}t^3\bigr)^2\bigl(1-{\rm e}^{w(\omega_2)}t\bigr)\bigl(1-{\rm e}^{w(\omega_2)}t^2\bigr)}\nonumber \\
\qquad=
\frac{w\bigl(d^{(2)}_{0}\bigr)}{1-t^3}
+\frac{w\bigl(d^{(2)}_{\omega_1}\bigr)}{1-{\rm e}^{w(\omega_1)}t^3}
+\frac{w\bigl(e^{(2)}_{\omega_1}\bigr)}{\bigl(1-{\rm e}^{w(\omega_1)}t^3\bigr)^2}
+\frac{w\bigl(d^{(2)}_{\omega_2}\bigr)}{1-{\rm e}^{w(\omega_2)}t}
+\frac{w\bigl(e^{(2)}_{\omega_2}\bigr)}{1-{\rm e}^{w(\omega_2)}t^2},\label{eqn:PsumG2pfd}
\end{gather}
where \smash{$d^{(2)}_{0},d^{(2)}_{\omega_1},e^{(2)}_{\omega_1},d^{(2)}_{\omega_2},e^{(2)}_{\omega_2}\in \K[t]$} are polynomials in $t$ of degrees at most 2, 2, 2, 0, 1, respectively, and are determined by the decomposition in the case $w=1$.

\begin{Proposition}\label{prop:pfdG2}
We have
\[
\mathcal{P}^{(2)}(t) =
\frac{1}{\Delta}\sum_{w\in W}(-1)^{\ell(w)}{\rm e}^{w(\rho)}
\left(
\frac{w\bigl(d^{(2)}_{0}\bigr)}{1-t^3}
+\frac{w\bigl(d^{(2)}_{\omega_1}\bigr)}{1-{\rm e}^{w(\omega_1)}t^3}
+\frac{w\bigl(d^{(2)}_{\omega_2}\bigr)}{1-{\rm e}^{w(\omega_2)}t}
\right),
\]
that is, the terms involving $w\bigl(e^{(2)}_{\omega_1}\bigr)$ and $w\bigl(e^{(2)}_{\omega_2}\bigr)$ cancel out after the Weyl group summation in~\eqref{eqn:PsumG2}.
\end{Proposition}

\begin{proof}
It is sufficient to show that for $\la\in \{\omega_1,\omega_2\}$,
\smash{$
\sum_{w\in W_{\la}}(-1)^{\ell(w)}{\rm e}^{w(\rho)}w\bigl(e^{(2)}_{\la}\bigr)=0$}.
The isotropy subgroup $W_{\la}$ has only two elements for each $\la$.
It is verified in the accompanying file using explicit expressions for the partial fraction coefficients in~\eqref{eqn:PsumG2pfd}.
\end{proof}
This completes step \ref{algo:resloc} of Section~\ref{sec:framework} for the $G_2$, $a=2$ case.

For $\mu\in \{0,\omega_1,\omega_2\}$ and $j\in \{0,1,2\}$, define the coefficients \smash{$a^{(2)}_{\mu;j}\in \K$} by
\[
	d^{(2)}_{0}=\sum_{j=0}^{2}a^{(2)}_{0;j}t^j,\qquad
	d^{(2)}_{\omega_1}=\sum_{j=0}^{2}a^{(2)}_{\omega_1;j}t^j,\qquad
	d^{(2)}_{\omega_2}=a^{(2)}_{\omega_2;0},
\]
and \smash{$a^{(2)}_{\omega_2;j}=0$} for $j=1,2$.
For $m\in \Z_{\geq 0}$, the coefficients of \smash{$\mathcal{P}^{(2)}(t)=\sum_{m=0}^{\infty}\PP2m t^m$} satisfy
\begin{gather*}
\PP{2}{3m} =
\frac{1}{\Delta}\sum_{w\in W}(-1)^{\ell(w)}{\rm e}^{w(\rho)}
\bigl(
	w\bigl(a^{(2)}_{0;0}\bigr)
	+w\bigl(a^{(2)}_{\omega_1;0}\bigr){\rm e}^{m w(\omega_1)}
	+w\bigl(a^{(2)}_{\omega_2;0}\bigr){\rm e}^{3m w(\omega_2)}
\bigr),\\
\PP{2}{3m+1} =
\frac{1}{\Delta}\sum_{w\in W}(-1)^{\ell(w)}{\rm e}^{w(\rho)}
\bigl(
	w\bigl(a^{(2)}_{0;1}\bigr)
	+w\bigl(a^{(2)}_{\omega_1;1}\bigr){\rm e}^{m w(\omega_1)}
	+w\bigl(a^{(2)}_{\omega_2;0}\bigr)e^{(3m+1) w(\omega_2)}
\bigr),\\
\PP{2}{3m+2} =
\frac{1}{\Delta}\sum_{w\in W}(-1)^{\ell(w)}{\rm e}^{w(\rho)}
\bigl(
	w\bigl(a^{(2)}_{0;2}\bigr)
	+w\bigl(a^{(2)}_{\omega_1;2}\bigr){\rm e}^{m w(\omega_1)}
	+w\bigl(a^{(2)}_{\omega_2;0}\bigr)e^{(3m+2) w(\omega_2)}
\bigr).
\end{gather*}
This follows by taking coefficients of \smash{$t^{3m}$}, \smash{$t^{3m+1}$}, and \smash{$t^{3m+2}$} in Proposition~\ref{prop:pfdG2} and expanding each denominator as a~geometric series.
It can be rewritten in the form
\begin{equation}\label{eqn:pp2}
\PP{2}{3m+j}=\sum_{\lambda\in \Lambda_2}D^{(2)}_{\lambda}{\rm e}^{j\lambda}{\rm e}^{3m\lambda}
\;+\sum_{\lambda\in \Lambda'_2}E^{(2)}_{\lambda,j}{\rm e}^{m\lambda},
\qquad j\in \{0,1,2\}.
\end{equation}
Here the coefficients $D^{(2)}_{\lambda}$ and $E^{(2)}_{\lambda,j}$ are given by
\[
D^{(2)}_{\lambda}:=
\frac{1}{\Delta}\sum_{\substack{w\in W\\ w(\omega_2)=\lambda}}
(-1)^{\ell(w)}{\rm e}^{w(\rho)}w\bigl(a^{(2)}_{\omega_2;0}\bigr),
\qquad \lambda\in \Lambda_2,
\]
\[
E^{(2)}_{\lambda,j}:=
\frac{1}{\Delta}\sum_{\substack{w\in W\\ w(\omega_1)=\lambda}}
(-1)^{\ell(w)}{\rm e}^{w(\rho)}w\bigl(a^{(2)}_{\omega_1;j}\bigr),
\qquad \lambda\in W\omega_1,\quad j\in \{0,1,2\},
\]
as well as
\[
E^{(2)}_{0,j}:=
\frac{1}{\Delta}\sum_{w\in W}
(-1)^{\ell(w)}{\rm e}^{w(\rho)}w\bigl(a^{(2)}_{0;j}\bigr),
\qquad j\in \{0,1,2\}.
\]

Let us carry out step~\ref{algo:rescompe} in a~slightly different way from that in Section~\ref{sec:f4}.
Instead of working with $C\bigl(\Pa,\la,\zeta,l\bigr)$ in
\smash{$
\Pam =\sum_{(\la,\zeta,l)} C\bigl(\Pa,\la,\zeta,l\bigr)\zeta^m {\rm e}^{m \la/l}$},
we just work with \smash{$D^{(1)}_{\lambda}$}, \smash{$D^{(2)}_{\lambda}$} and~\smash{$E^{(2)}_{\lambda,j}$} in~\eqref{eqn:pp1} and~\eqref{eqn:pp2} to avoid the use of root of unity in the actual computation, and work with rational functions in $x_i={\rm e}^{\alpha_i}$ with coefficients in $\QQ$.
We first verify that the family~\smash{$\bigl(\Pam\bigr)$} satisfies the $Q$-system by checking the relations among the coefficients appearing above; this verification is carried out in the accompanying file.
Once this is established, it suffices to check the equality \smash{$\Pam=\Qam$} only for $m=0$ and $m=1$.
Since this condition is satisfied in the present case (that is, \smash{$\Pam$} gives the correct decomposition of the KR modules), it follows that \smash{$\Pam=\Qam$} holds for all $m\in\Z_{\geq 0}$.
This is because \smash{$\bigl(\Qam\bigr)$} is a~solution of the $Q$-system and~\smash{$\Qam\neq 0$}, and hence is uniquely determined by the values at $m=0,1$.
This finishes the proof of Theorem~\ref{thm:main2}.

\subsection*{Acknowledgements}

The author thanks the anonymous referees for their careful reading and helpful comments.
The author is supported by a~KIAS Individual Grant (SP067302) via the June E Huh Center for Mathematical Challenges at Korea Institute for Advanced Study.

\pdfbookmark[1]{References}{ref}
\LastPageEnding

\end{document}